\documentclass[12pt]{amsart}
\usepackage[margin=1in]{geometry}
\usepackage{amsmath,amsfonts,amsthm,amssymb,bbm}
\usepackage{graphicx,color,dsfont}
\usepackage{enumitem}
\usepackage{fourier}

\begin{document}

\newtheorem{theorem}{Theorem}
\newtheorem{lemma}{Lemma}
\newtheorem{proposition}{Proposition}
\newtheorem{rmk}{Remark}
\newtheorem{example}{Example}
\newtheorem{exercise}{Exercise}
\newtheorem{definition}{Definition}
\newtheorem{corollary}{Corollary}
\newtheorem{notation}{Notation}
\newtheorem{claim}{Claim}

\newtheorem{dif}{Definition}

 \newtheorem{thm}{Theorem}[section]
 \newtheorem{cor}[thm]{Corollary}
 \newtheorem{lem}[thm]{Lemma}
 \newtheorem{prop}[thm]{Proposition}
 \theoremstyle{definition}
 \newtheorem{defn}[thm]{Definition}
 \theoremstyle{remark}
 \newtheorem{rem}[thm]{Remark}
 \newtheorem*{ex}{Example}
 \numberwithin{equation}{section}

\newcommand{\vertiii}[1]{{\left\vert\kern-0.25ex\left\vert\kern-0.25ex\left\vert #1
    \right\vert\kern-0.25ex\right\vert\kern-0.25ex\right\vert}}

\newcommand{\R}{{\mathbb R}}
\newcommand{\C}{{\mathbb C}}
\newcommand{\U}{{\mathcal U}}
\newcommand{\norm}[1]{\left\|#1\right\|}
\renewcommand{\(}{\left(}
\renewcommand{\)}{\right)}
\renewcommand{\[}{\left[}
\renewcommand{\]}{\right]}
\newcommand{\f}[2]{\frac{#1}{#2}}
\newcommand{\im}{i}
\newcommand{\cl}{{\mathcal L}}
\newcommand{\ck}{{\mathcal K}}

\newcommand{\al}{\alpha}
\newcommand{\vro}{\varrho}
\newcommand{\be}{\beta}
\newcommand{\wh}[1]{\widehat{#1}}
\newcommand{\ga}{\gamma}
\newcommand{\Ga}{\Gamma}
\newcommand{\de}{\delta}
\newcommand{\ben}{\beta_n}
\newcommand{\De}{\Delta}
\newcommand{\ve}{\varepsilon}
\newcommand{\ze}{\zeta}
\newcommand{\Th}{\Theta}
\newcommand{\ka}{\kappa}
\newcommand{\la}{\lambda}
\newcommand{\laj}{\lambda_j}
\newcommand{\lak}{\lambda_k}
\newcommand{\La}{\Lambda}
\newcommand{\si}{\sigma}
\newcommand{\Si}{\Sigma}
\newcommand{\vp}{\varphi}
\newcommand{\om}{\omega}
\newcommand{\Om}{\Omega}
\newcommand{\ra}{\rightarrow}

\newcommand{\ro}{{\mathbf R}}
\newcommand{\rn}{{\mathbf R}^n}
\newcommand{\rd}{{\mathbf R}^d}
\newcommand{\rmm}{{\mathbf R}^m}
\newcommand{\rone}{\mathbb R}
\newcommand{\rtwo}{\mathbf R^2}
\newcommand{\rthree}{\mathbf R^3}
\newcommand{\rfour}{\mathbf R^4}
\newcommand{\ronen}{{\mathbf R}^{n+1}}
\newcommand{\ku}{\mathbf u}
\newcommand{\kw}{\mathbf w}
\newcommand{\kf}{\mathbf f}
\newcommand{\kz}{\mathbf z}

\newcommand{\N}{\mathbf N}

\newcommand{\tn}{\mathbf T^n}
\newcommand{\tone}{\mathbf T^1}
\newcommand{\ttwo}{\mathbf T^2}
\newcommand{\tthree}{\mathbf T^3}
\newcommand{\tfour}{\mathbf T^4}

\newcommand{\zn}{\mathbf Z^n}
\newcommand{\zp}{\mathbf Z^+}
\newcommand{\zone}{\mathbf Z^1}
\newcommand{\zz}{\mathbf Z}
\newcommand{\ztwo}{\mathbf Z^2}
\newcommand{\zthree}{\mathbf Z^3}
\newcommand{\zfour}{\mathbf Z^4}

\newcommand{\hn}{\mathbf H^n}
\newcommand{\hone}{\mathbf H^1}
\newcommand{\htwo}{\mathbf H^2}
\newcommand{\hthree}{\mathbf H^3}
\newcommand{\hfour}{\mathbf H^4}

\newcommand{\cone}{\mathbf C^1}
\newcommand{\ctwo}{\mathbf C^2}
\newcommand{\cthree}{\mathbf C^3}
\newcommand{\cfour}{\mathbf C^4}
\newcommand{\dpr}[2]{\langle #1,#2 \rangle}

\newcommand{\sn}{\mathbf S^{n-1}}
\newcommand{\sone}{\mathbf S^1}
\newcommand{\stwo}{\mathbf S^2}
\newcommand{\sthree}{\mathbf S^3}
\newcommand{\sfour}{\mathbf S^4}

\newcommand{\lp}{L^{p}}
\newcommand{\lppr}{L^{p'}}
\newcommand{\lqq}{L^{q}}
\newcommand{\lr}{L^{r}}
\newcommand{\echi}{(1-\chi(x/M))}
\newcommand{\chip}{\chi'(x/M)}

\newcommand{\wlp}{L^{p,\infty}}
\newcommand{\wlq}{L^{q,\infty}}
\newcommand{\wlr}{L^{r,\infty}}
\newcommand{\wlo}{L^{1,\infty}}

\newcommand{\lprn}{L^{p}(\rn)}
\newcommand{\lptn}{L^{p}(\tn)}
\newcommand{\lpzn}{L^{p}(\zn)}
\newcommand{\lpcn}{L^{p}(\cn)}
\newcommand{\lphn}{L^{p}(\cn)}

\newcommand{\lprone}{L^{p}(\rone)}
\newcommand{\lptone}{L^{p}(\tone)}
\newcommand{\lpzone}{L^{p}(\zone)}
\newcommand{\lpcone}{L^{p}(\cone)}
\newcommand{\lphone}{L^{p}(\hone)}

\newcommand{\lqrn}{L^{q}(\rn)}
\newcommand{\lqtn}{L^{q}(\tn)}
\newcommand{\lqzn}{L^{q}(\zn)}
\newcommand{\lqcn}{L^{q}(\cn)}
\newcommand{\lqhn}{L^{q}(\hn)}

\newcommand{\lo}{L^{1}}
\newcommand{\lt}{L^{2}}
\newcommand{\li}{L^{\infty}}
\newcommand{\beqn}{\begin{eqnarray*}}
\newcommand{\eeqn}{\end{eqnarray*}}
\newcommand{\pplus}{P_{Ker[\cl_+]^\perp}}

\newcommand{\co}{C^{1}}
\newcommand{\ci}{C^{\infty}}
\newcommand{\coi}{C_0^{\infty}}

\newcommand{\ca}{\mathcal A}
\newcommand{\cs}{\mathcal S}
\newcommand{\cm}{\mathcal M}
\newcommand{\cf}{\mathcal F}
\newcommand{\cb}{\mathcal B}
\newcommand{\ce}{\mathcal E}
\newcommand{\cd}{\mathcal D}
\newcommand{\cn}{\mathcal N}
\newcommand{\cz}{\mathcal Z}
\newcommand{\crr}{\mathbf R}
\newcommand{\cc}{\mathcal C}
\newcommand{\ch}{\mathcal H}
\newcommand{\cq}{\mathcal Q}
\newcommand{\cp}{\mathcal P}
\newcommand{\cx}{\mathcal X}
\newcommand{\eps}{\epsilon}

\newcommand{\pv}{\textup{p.v.}\,}
\newcommand{\loc}{\textup{loc}}
\newcommand{\intl}{\int\limits}
\newcommand{\iintl}{\iint\limits}
\newcommand{\dint}{\displaystyle\int}
\newcommand{\diint}{\displaystyle\iint}
\newcommand{\dintl}{\displaystyle\intl}
\newcommand{\diintl}{\displaystyle\iintl}
\newcommand{\liml}{\lim\limits}
\newcommand{\suml}{\sum\limits}
\newcommand{\ltwo}{L^{2}}
\newcommand{\supl}{\sup\limits}
\newcommand{\df}{\displaystyle\frac}
\newcommand{\p}{\partial}
\newcommand{\Ar}{\textup{Arg}}
\newcommand{\abssigk}{\widehat{|\si_k|}}
\newcommand{\ed}{(1-\p_x^2)^{-1}}
\newcommand{\tT}{\tilde{T}}
\newcommand{\tV}{\tilde{V}}
\newcommand{\wt}{\widetilde}
\newcommand{\Qvi}{Q_{\nu,i}}
\newcommand{\sjv}{a_{j,\nu}}
\newcommand{\sj}{a_j}
\newcommand{\pvs}{P_\nu^s}
\newcommand{\pva}{P_1^s}
\newcommand{\cjk}{c_{j,k}^{m,s}}
\newcommand{\Bjsnu}{B_{j-s,\nu}}
\newcommand{\Bjs}{B_{j-s}}
\newcommand{\Ly}{\cl_+i^y}
\newcommand{\dd}[1]{\f{\partial}{\partial #1}}
\newcommand{\czz}{Calder\'on-Zygmund}
\newcommand{\chh}{\mathcal H}

\newcommand{\lbl}{\label}
\newcommand{\beq}{\begin{equation}}
\newcommand{\eeq}{\end{equation}}
\newcommand{\beqna}{\begin{eqnarray*}}
\newcommand{\eeqna}{\end{eqnarray*}}
\newcommand{\bp}{\begin{proof}}
\newcommand{\ep}{\end{proof}}
\newcommand{\bprop}{\begin{proposition}}
\newcommand{\eprop}{\end{proposition}}
\newcommand{\bt}{\begin{theorem}}
\newcommand{\et}{\end{theorem}}
\newcommand{\bex}{\begin{Example}}
\newcommand{\eex}{\end{Example}}
\newcommand{\bc}{\begin{corollary}}
\newcommand{\ec}{\end{corollary}}
\newcommand{\bcl}{\begin{claim}}
\newcommand{\ecl}{\end{claim}}
\newcommand{\bl}{\begin{lemma}}
\newcommand{\el}{\end{lemma}}
\newcommand{\dea}{(-\De)^\be}
\newcommand{\naa}{|\nabla|^\be}
\newcommand{\cj}{{\mathcal J}}
\newcommand{\ubb}{{\mathbf u}}

\title[Waves in the NLS system of  the third-harmonic generation]{On the stability of solitary waves in the NLS system of  the third-harmonic generation}

\author[Abba Ramadan]{\sc Abba Ramadan}
\address{ Department of Mathematics,
The University of Alabama, Room 150, Gordon Palmer Hall,
 505 Hackberry Lane,
Tuscaloosa, AL 35401, USA}
\email{aramadan@ua.edu}

\author[Atanas G. Stefanov]{\sc Atanas G. Stefanov}
\address{ Department of Mathematics,
	University of Alabama - Birmingham, 
	University Hall, Room 4005, 
	1402 10th Avenue South
	Birmingham AL 35294-1241
	 }
\email{stefanov@uab.edu}



\date{\today}
 
\begin{abstract}
 We consider the NLS system of  the third-harmonic generation, which was  introduced in \cite{Sammut}. Our interest is in solitary wave solutions and their stability properties. The recent work of Oliveira and Pastor, \cite{OP}, discussed   global well-posedness vs. finite time blow up, as well as  other aspects of the dynamics. These  authors have also constructed solitary wave solutions, via the method of mountain pass/Nehari manifold,  in an appropriate range of parameters. Specifically, the  waves exist only in spatial dimensions $n=1,2,3$. They have also establish some stability/instability results for these waves. 
 
 In this work, we systematically build and study solitary waves for  this important  model. We construct the waves in the largest possible parameter space, and we provide a complete classification of their stability.  In dimension one, we show stability, whereas in $n=2,3$, they are generally spectrally unstable, except for a small region, where they do enjoy an extra pseudo-conformal symmetry. 
 
 Finally,  we discuss instability by blow-up. In the case $n=3$, and for more restrictive set of parameters,  we use virial identities methods to derive the strong instability, in the spirit of Ohta's approach, \cite{O}. In $n=2$,  the virial identities reduce matters, via conservation of mass and energy, to the initial data.  Our conclusions mirror closely the well-known results for the scalar cubic  focussing NLS, while the proofs are  much more involved.

\end{abstract}

\thanks{ Stefanov is partially supported by NSF-DMS \# 1908626 and NSF-DMS \# 2204788.} 

\maketitle

\section{Introduction}
\label{sec:1}
In this paper, we are interested in the following system of non-linear Schr\"odinger equations 
  \begin{equation}
  	\label{10}
  	\begin{cases} 
  		iu_t+\De u-u+(\frac{1}{9}|u|^2+2|v|^2)u+\frac{1}{3}\bar{u}^2 v=0 \\
  		i\sigma  v_t+\De v-\mu v+(9|v|^2+2|u|^2)v +\frac{1}{9}u^3=0,
  	\end{cases}
  \end{equation}
where $\mu, \si>0$. 
  This system arises in nonlinear optics applications. Specifically,  this model was proposed in \cite{Sammut}, where the interaction of a 
  a monochromatic beam  with frequency $\om$ propagates in a cubic (Kerr) medium and interacts with its third harmonic, frequency $3\om$. In such a scenario, consider the electric field $E$, which satisfies the Maxwell's equation. 
  
  Denoting $U, V$ the amplitudes of the   modes, corresponding to frequencies $\om, 3\om$, a system very similar to \eqref{10} arises as an envelope approximation of the Maxwell's equation. After some rescaling and non-dimensionalizaation, one obtains  the model \eqref{10}. For more details, the reader should consult \cite{Sammut} and also \cite{OP}. 
  
  A few other noteworthy details about the system \eqref{10} are as follows.  It is a Hamiltonian system, and as such is generated by Hamiltonian function, and, at least formally,  conserves energy and mass. The conserved quantities take the form 
  \beqn
  E(u,v) &=& \frac{1}{2}\int_{\rn} (|\nabla u|^2+
  |\nabla v|^2+|u|^2+\mu|v|^2)dx-\int_{\rn}  \frac{1}{36}|u|^4+\frac{9}{4}|v|^4+|u|^2|v|^2+\frac{1}{9}
  \Re(\bar{u}^3 v)dx\\
  M(u,v) &=& \int_{\rn}  (|u|^2+3\sigma|v|^2)dx
  \eeqn
  where  $u, v$ are classical solutions of \eqref{10}.

  Our goal is to construct  various solitary waves for the system \eqref{10} and to study their stability properties. To this end, we are looking for solutions in the ansatz  $u=e^{ i \om t} P, v=e^{ 3 i \om t} Q$. Note that this form of time dependence directly corresponds to the model, which seeks to model the third harmonic generation. Mathematically, this intrinsic property  manifests itself as this specific  symmetry of the model, which factors out the time dependence and naturally leads to the  following  elliptic system 
    \begin{equation}
   	\label{20}   	
   	\begin{cases} 
   		-\De P+(\om+1)P -\left[\frac{1}{9}P^3+2Q^2 P+\frac{1}{3}P^2Q\right]=0 \\
   		-\De Q+(\mu+3\si \om) Q -\left[9Q^3+2P^2Q+\frac{1}{9}P^3\right]=0.
   	\end{cases}
   \end{equation}
 First we state the following, so-called  Pohozaev identities, already derived in the earlier paper  \cite{OP}.

\begin{lemma}(Pohozhaev's identities)
	\label{Poh}
	
	Assume that \eqref{20} has solution $(P,Q)\in H^1(\rn)\times H^1(\rn)).$ Then the following identities hold:
	\begin{eqnarray}
		\label{Poh1}
	& & 	\int\left(|\nabla P|^2+(\om+1)P^2\right)dx=\int \left(\frac{1}{9}P^4+2P^2Q^2+\frac{1}{3}P^3Q\right)dx,  \\ 
	& & 	\label{Poh2}
		\int(|\nabla Q|^2+(\mu+3\sigma\om)Q^2)dx=\int\left( 9Q^4+2P^2Q^2+\frac{1}{9}P^3Q\right)dx.
	\end{eqnarray}
In addition, 
	\begin{equation}
		\label{Poh3}
		(4-n)\int(|\nabla Q|^2+|\nabla P|^2)dx=
		n(\om+1)\int P^2 dx+n(\mu+3\sigma\om)\int Q^2dx. 
	\end{equation}
	
\end{lemma}   
Obverve that for $n\geq 4$ and 
$\om> \max\{-1,-\frac{\mu}{3\sigma}\}$ the solitary waves for the system do not exist in $(H^1(\rn)\times H^1(\rn))$.  
Next, we discuss the linearization of \eqref{10} around the standing waves 
$u=e^{i \om t}P,v=e^{3 i \om t}Q$. 
\subsection{The linearized system}
We perform a standard linearization procedure, namely we take 
$u=e^{ i \omega t}[P+\phi],v=e^{ 3 i \omega t}[Q+\psi]$, plugging it in \eqref{10} and ignoring the higher order terms $O(\phi^2,\psi^2,\psi\phi)$, we arrive at the linearized system, which after $\phi =(\Re \phi, \Im \phi) =: (\phi_1, \phi_2)$ and  $\psi =(\Re \psi, \Im \psi) =: (\psi_1, \psi_2)$   can be written as  follows 
\begin{equation}
\label{412} 
\left(\begin{array}{c}
\ \phi_1\\ \psi_1 \\ \phi_2\\ \psi_2
\end{array}\right)_t = \left(\begin{array}{cccc}
0 & 0 & 1 & 0 \\ 0 & 0 & 0 & \frac{1}{\sigma}\\ -1 & 0 & 0 & 0\\0 & -\frac{1}{\sigma} & 0 & 0
\end{array}\right) \left(\begin{array}{cccc}
\cl_1 & A & 0 & 0 \\ A & \cl_3 & 0 & 0\\ 0 & 0 & \cl_2 & B \\ 0 & 0 & B & \cl_4
\end{array}\right) \left(\begin{array}{c}
 \phi_1\\ \psi_1 \\ \phi_2\\ \psi_2
\end{array}\right).
\end{equation}
Here, we have introduced  the following scalar Schr\"odinger operators
\begin{eqnarray*}
\cl_1 &=&-\De +(\om+1) - \left(\frac{P^2}{3}+2Q^2+\frac{2}{3}PQ\right),  \ \ 
\cl_2 = -\De +(\om+1) -  \left(\frac{P^2}{9}+2Q^2+\frac{2}{3}PQ\right) ,\\
\cl_3 &=&-\De +(\mu+3\sigma\om) - (27Q^2+2P^2),  \ \ 
\cl_4 = -\De +(\mu+3\sigma\om) -  (9Q^2+2P^2) 
\end{eqnarray*}
and the potentials are $A=-4PQ-\frac{P^2}{3},B=-\frac{P^2}{3}$. 
Introduce the operators $\cl_+,\cl_-$as follows 
\begin{equation}
	\label{414} 
	\cl_+:=\begin{bmatrix}
		\cl_1 & A \\
		A & \cl_3
	\end{bmatrix}, 
	\ \ \cl_-:=\begin{bmatrix}
		\cl_2 & B \\
		B & \cl_4
	\end{bmatrix}.
\end{equation}
Note that, by direct inspection, we have $\cl_-\left(\begin{array}{c}
	P\\ 3 Q 
\end{array}\right)=0$. Also, by 
taking the spatial derivative of \eqref{20} we have
\begin{eqnarray}
	\label{a:19} 
	\cl_+\left(\begin{array}{c}
		P' \\ Q'
	\end{array}\right)=0. 
\end{eqnarray} 
 This is not surprising due to translation and modulational  invariance of \eqref{10}.
Also assigning $\begin{pmatrix}
\phi_1 \\ \psi_1 \\ \phi_2 \\\psi_2
\end{pmatrix}\to e^{\lambda t}\begin{pmatrix}
v_1 \\v_2 \\ v_3 \\v_4
\end{pmatrix}=:e^{\lambda t}\vec{v}$, 
we obtain the following time-independent   eigenvalue problem
\begin{equation}
\label{Eigen}
\mathcal{J}\cl\vec{v}=\lambda\vec{v}.
\end{equation}
Here 
$$
\mathcal{J}:= \left(\begin{array}{cccc}
0 & 0 & 1 & 0 \\ 0 & 0 & 0 & \frac{1}{\sigma}\\ -1 & 0 & 0 & 0\\0 & -\frac{1}{\sigma} & 0 & 0
\end{array}\right), \ \ \cl:=\begin{pmatrix}
\cl_+ & \textbf{0} \\
\textbf{0} & \cl_-
\end{pmatrix}.
$$
In a standard way, we say that the time-periodic waves $(e^{i \om t} P, e^{3 i \om t} Q)$ are {\it spectrally stable} if the eigenvalue problem \eqref{Eigen} does not have solutions $(\la, \vec{v}): \Re\la>0, \vec{v}\in Dom(\cl)=(H^2(\rn))^4$. 

\subsection{Main results: existence}
Our first group of results deals with the existence of appropriate solutions of \eqref{20}. We will provide at least two ways of producing solutions of \eqref{20}, both variational in nature.
\begin{theorem}
	\label{theo:10}
Let $1\leq n\leq 3$, $\si>0$, and $\om>\max(-1, -\f{\mu}{3\si})$. Then, there exist bell-shaped exponentally decaying solutions $(P,Q)\in (H^1(\rn)\times H^1(\rn))$ of \eqref{20}. 
In addition, the self-adjoint linearized operators $\cl_\pm$ satisfy the spectral properties: 
$\cl_-\geq 0$, while $\cl_+$ has exactly one negative eigenvalue. 
\end{theorem}
{\bf Remark:} Note that the condition $\om>\max(-1, -\f{\mu}{3\si})$, and $1\leq n\leq 3$ is necessary for the existence of localized pair of functions $P,Q$ satisfying \eqref{20}. 

Next, we provide an additional construction of solutions $(P,Q)$ of \eqref{20}, the so-called normalized waves. These only exist  in the case $n=1$. These solutions are somewhat special in several ways to be described in the sequel. 
\begin{theorem}
	\label{theo:40} 
	Let $n=1$, $\si>0$. Then, for each $\la>0$, there exists $\om=\om_\la>: \om>\max(-1, -\f{\mu}{3\si})$ and a pair of bell-shaped functions $(P,Q): \int_{\rone} (P^2(x) +3\si Q^2(x) ) dx =\la$, so that \eqref{20} holds true. Moreover, $\cl_-\geq 0$, while $n(\cl_+)=1$. In fact, 
	$$
		\dpr{\cl_+h}{h}\geq 0, \forall h \perp  \left(\begin{matrix}
				P \\ 3 \si Q
			\end{matrix}\right). 
	$$
The corresponding waves $(e^{i \om_\la t} P, e^{3 i \om_\la t} Q)$ are spectrally stable.  
\end{theorem}
Our next result is about the stability of generic solutions of \eqref{20}. 
\subsection{Main results: Spectral Stability} 
We start with the higher dimensional case $n>1$, which mostly features unstable solitons. 
\begin{theorem}
	\label{theo:20} 
	Let $\si>0$, $n=2,3$, $\om>\max(-1, -\f{\mu}{3\si})$.  Assume that $(P,Q)\in (H^1(\rn)\times H^1(\rn))$ satisfies  \eqref{20}, subject to the spectral property: $\cl_-\geq 0$,  $\cl_+$ has exactly one negative eigenvalue. 
	
For $n=3$, the solution $(e^{i \om t} P, e^{3 i \om t} Q)$ is spectrally unstable. 
For $n=2$, $\mu\neq 3\si$, the solution $(e^{i \om t} P, e^{3 i \om t} Q)$ is spectrally unstable.

For $n=2, \mu=3\si$,  the solution $(e^{i \om t} P, e^{3 i \om t} Q)$ is in fact spectrally stable.
	
\end{theorem}
{\bf Remarks:} 
\begin{itemize}
	\item The property $\cl_-\geq 0$ follows  easily from the Perron-Frobenius property of $\cl_-$, once we assume $P>0, Q>0$. In fact, this condition is also necessary for $\cl_-\geq 0$ . 
	\item The requirement $n(\cl_+)=1$ is satisfied for the variational solutions constructed in Theorem \ref{theo:10}, but they do not have to be constructed variationally or under some specific variational procedure. 
	\item The spectral stability of the solution in the case $n=2, \mu=3\si$ is via a mechanism similar to the stability of soliton the 2D cubic NLS. Specifically, the algebraic multiplicity of  the zero eigenvalue of the associated e-value problem is exactly eight, vis-a-vis the generic value\footnote{For any $p$ power nonlinearity $p\neq 3$}  of six.   
	
\end{itemize}
\subsection{Main results: Instability by blow-up}
Next, we give definition of strong instability, also known as instability by 
blow-up. 
 \begin{definition}
	\label{22}
	We say $(e^{i \om t} P, e^{3 i \om t} Q)$ is strongly unstable if for any $\epsilon>0$ there exist $(u_0,v_0)\in (H^1(\rn)\times H^1(\rn)) $ such that $\|(u_0,v_0)-(P,Q)\|_{(H^1(\rn)\times H^1(\rn)) }<\epsilon$ and the solution $(u(t),v(t))$ of \eqref{10} with $(u(0),v(0))=(u_0,v_0)$ blows up in finite time.
\end{definition}
 
\begin{theorem}
    \label{21}
    Let $\si=3$, $\mu=3\si$, and $n=2,3$, $\om>\max(-1, -\f{\mu}{3\si})$. Assume that $(P,Q)\in (H^1(\rn)\times H^1(\rn))$ satisfies  \eqref{20}, but in addition, they are minimizers of a constrained variational problem, see Proposition \ref{prop:22} below. Then, the solitary wave solution $(e^{i \om t} P, e^{3 i \om t} Q)$ of \eqref{10} is strongly unstable. 
\end{theorem}
 We now plan ahead for our presentation as follows. In Section \ref{sec:2}, we construct the solitary waves via the Weinstein's functional approach, in all the relevant dimensions $n=1,2,3$. We also show some relevant spectral properties for the linearized operators $\cl_\pm$. In addition, the normalized waves are constructed in dimension $n=1$ only, as this is the only case where this is possible. In particular, we show that the waves are spectrally stable. In Section \ref{sec:3}, we  rigorously establish the instability of the waves in $n=2,3$. Finally, in Section \ref{sec:4}, we show, by means of virial identities in the spirit of Ohta, that the spectrally unstable waves for the two and three dimensional case, are in fact strongly unstable.

   \section{Existence and spectral properties of the waves}
   \label{sec:2} 
   We start with some preliminaries. 
   \subsection{Preliminaries} We work with standard function spaces like $L^p, W^{s,p}$ etc. 
   Since a particular interest will be placed on bell-shaped solutions, we go ahead with the standard introductions. 
   
   First, it is well-known that 
   or each $L^1_{loc}(\rn)$ function $f$, there is a decreasing rearrangement $f^*$. More precisely, there is unique function $f^*\in L^1_{loc}(\rn)$, so that for each $\la>0$, so that $f, f^*$ are equi-measurable ( that is, for each $\la>0$, 
   $|\{x: |f(x)|>\la\}|=|\{x: |f^*(x)|>\la\}|$) and so that $f^*(x)=q(|x|)$, where $q:\rone_+\to \rone_+$ is a non-increasing function. Next, we introduce the notion of bell-shaped function. 
   \begin{definition}
   	\label{defi:10} 
   	We say that a function $f\in L^1_{loc}(\rn)$ is bell-shaped, if $f=f^*$. 
   \end{definition}
   It is standard that for $p: 1\leq p\leq \infty$, $\|f\|_{L^p}=\|f^*\|_{L^p}$. Also, we record the standard rearrangement and Szeg\"o  inequalities 
   \begin{eqnarray}
   	\label{40} 
   	& & 	\int_{\rn} f g dx \leq 	\int_{\rn} f^* g^* dx \\
   	\label{50} 
   	& & \|\nabla f\|_{L^2(\rn)}\geq \|\nabla f^*\|_{L^2(\rn)}.
   \end{eqnarray}
  
   \subsection{Existence of the waves via Weinstein's functionals}
   Next,  we employ a variational method to produce waves. The goal is to find an effective way to construct waves that satisfy \eqref{20}, after some straightforward transformation such as scaling etc. 
   As we would like to cover an existence theory for a wider range of parameters, we fix $\al, \be>0$ and consider the following Weinstein functional 
   \begin{equation}
   	\label{42} 
   	W_{\al,\be}[u,v]:=\f{\int\limits_{\rn} (|\nabla u|^2+
   		|\nabla v|^2+\al |u|^2+\be |v|^2)dx}{\left(\int\limits_{\rn} \frac{1}{36}|u|^4+\frac{9}{4}|v|^4+|u|^2|v|^2+\frac{1}{9}
   		\Re(\bar{u}^3 vdx)\right)^{\f{1}{2}} }\to \min
   \end{equation}
where the minimum is taken over $u,v\in \cs: (u,v)\neq (0,0)$.  Note the homogeneity of the functional 
$W_{\al, \be}$, in the sense that for each $a\neq 0$, $W_{\al, \be}[a u, a v]=W_{\al, \be}[u, v]$. Thus, it is clear that solving the uncosntrained variational problem \eqref{42} is equivalent to solving the following constrained maximization problem 
   \begin{equation}
   	\label{30}
   	\begin{cases}
   		J[u,v]:=\int\limits_{\rn}  \left(\frac{1}{36}|u|^4+\frac{9}{4}|v|^4+|u|^2|v|^2+\frac{1}{9}
   		\Re(\bar{u}^3 v)\right) dx \to \max \\
   	I[u,v]:=	\int\limits_{\rn} (|\nabla u|^2+
   		|\nabla v|^2+\al |u|^2+\be |v|^2)dx=1.
   	\end{cases}
   \end{equation}

Our main existence result is the following. 
\begin{proposition}
	\label{prop:10} 
	Let $\al, \be>0$ and $1\leq n\leq 3$. Then, the variational problem \eqref{30} has solution $(U,V)$, which is necessarily bell-shaped. In addition, it satisfies the Euler-Lagrange equations 
	\begin{equation}
   	\label{E-L}   	
   	\begin{cases} 
   		-\De U+\alpha U -C(\alpha,\beta)\left[\frac{1}{9}U^3+2V^2 U+\frac{1}{3}U^2 V\right]=0 \\
   		-\De V+\beta V -C(\alpha,\beta)\left[9V^3+2U^2 V+\frac{1}{9}U^3\right]=0,
   	\end{cases}
   \end{equation}
for some positive scalar $C(\al, \be)$. 
Finally, the  linearized operators
\begin{equation}
\label{cpl}
L_+:=\begin{bmatrix}
L_1 & C(\alpha,\beta)A \\
C(\alpha,\beta)A & L_3
\end{bmatrix}, L_-:=\begin{bmatrix}
L_2 & C(\alpha,\beta)B \\
C(\alpha,\beta)B & L_4
\end{bmatrix}
\end{equation}	
	where 
	 $$
	 L_1:=-\Delta+\alpha-C(\alpha,\beta)(\frac{U^2}{3}+2V^2+\frac{2}{3}U V), \ \ L_3:=-\Delta+\beta-C(\alpha,\beta)(U^2+27V^2)
	 $$
	 (and similarly $L_2, L_4$) satisfies that for each test function $h: h\perp \begin{bmatrix}
(-\De+\al) U \\
(-\De+\be) V
\end{bmatrix}$, then  $\langle\cl_+h,h\rangle\geq 0.$ 
\end{proposition}
\begin{proof}
	It is not hard to see that the variational problem \eqref{30} is equivalent, and has the same solutions if any, with the following one 
	 \begin{equation}
		\label{31}
		\begin{cases}
			\int\limits_{\rn}  \left(\frac{1}{36}|u|^4+\frac{9}{4}|v|^4+|u|^2|v|^2+\frac{1}{9}
			\Re(\bar{u}^3 v)\right) dx \to \max \\
			\int\limits_{\rn} (|\nabla u|^2+
			|\nabla v|^2+\al |u|^2+\be |v|^2)dx\leq 1.
		\end{cases}
	\end{equation}
In fact, we claim that if $(u,v)$ are functions which satisfy $\int\limits_{\rn} (|\nabla u|^2+
|\nabla v|^2+\al |u|^2+\be |v|^2)dx<1$, we can take $(\tilde{u}, \tilde{v}):=a (u,v)$, where $a=\left(\int\limits_{\rn} (|\nabla u|^2+
|\nabla v|^2+\al |u|^2+\be |v|^2)dx\right)^{-\f{1}{2}}>1$. Clearly $(\tilde{u}, \tilde{v})$ still satisfies the constraint 
$I[\tilde{u}, \tilde{v}]=1$, while $J[\tilde{u}, \tilde{v}]=a^4 J[u,v]$ and hence it is better choice than $(u,v)$. 

Next, by H\"older's and Sobolev embedding, 
we see that for $1\leq n\leq 3$, we have that 
$$
J[u,v]\leq C (\|u\|_{L^4}^4+\|v\|_{L^4}^4))\leq C\|(u,v)\|_{H^{\f{n}{4}}\times H^{\f{n}{4}}}^4\leq C\|(u, v)\|_{H^1\times H^1}^4\leq C_{\al, \be} (I[u,v])^2.
$$
This shows that the problem \eqref{42} (and equivalently \eqref{30}) is well-posed, i.e. 
$$
J_{\max}:= \sup_{I]u,v]=1} \int\limits_{\rn}  \left(\frac{1}{36}|u|^4+\frac{9}{4}|v|^4+|u|^2|v|^2+\frac{1}{9}
\Re(\bar{u}^3 v)\right) dx<\infty
$$
is well-defined. 
Pick a  maximizing sequences, that is $(u_k, v_k)$, so that 
$
J[u_k, v_k]\to J_{\max}.
$
It is now  clear that the corresponding decreasing rearrangements $u_k^*, v_k^*$ are a better alternative. 
Indeed, due to \eqref{50}, we have that $I[u_k^*, v_k^*]\leq I[u_k, v_k]=1$, while by \eqref{40}, we have $J[u_k^*, v_k^*]\geq J[u_k, v_k]$. It follows that $(u_k^*, v_k^*)$ is a maximizing sequence as well. We now show that this sequence  is compact in 
$L^4(\rn)\times L^4(\rn) $. Specifically, by the Riesz-Kolmogorov compactness criteria and since\footnote{and here is where we need $n\leq 3$, as $n=4$ does not guarantee compactness of the embedding anymore}    for $1\leq n\leq 3$, 
$H^1(\rn)\hookrightarrow   L^4(\rn)$, it suffices to check that 
\begin{equation}
	\label{62} 
	\lim_{R\to \infty} \sup_k \int_{|x|>R} |u_k^*|^4+ |v_k^*|^4 dx=0.
\end{equation}
This actually follows easily from the bell-shapedness, since for every $M>1$, and $|x|=M$, we have 
$$
|u_k^*(x)|^2 \leq \f{1}{|B_M|} \int_{B_M} |u_k^*(y)|^2 dy \leq C M^{-n} \|u_k^*\|_{L^2}^2\leq C |x|^{-n},
$$
as $ \|u_k^*\|_{L^2}^2\leq C I[u_k^*, v_k^*]=C$. 
Similar estimate holds for $v_k^*$. It follows that 
$$
\int_{|x|>R} |u_k^*|^4+ |v_k^*|^4 dx\leq C \int_{|x|>R} |x|^{-2 n} dx = C R^{-n},
$$
whence \eqref{62} follows.

Next is to derive the Euler Lagrange equation \eqref{E-L}. To that end take any arbitary test function $h=[h_1,h_2], \epsilon>0$ we have 
\begin{eqnarray*}
& & J(U+\epsilon h_1,V+\epsilon h_2) =  J(U,V)+\epsilon\left(\langle \frac{1}{9}U^3+2U V^2+\frac{1}{3}U^2V,h_1\rangle)+\langle 9V^3+2U^2V+\frac{U^3}{9},h_2\rangle \right) +\\
&+& \f{\epsilon^2}{2}(\langle \left(\frac{U^2}{3}+2V^2 +\frac{2}{3}U V \right)h_1,h_1\rangle 
+\dpr{(8U V+\frac{2}{3}U^2)h_1}{h_2}+ \langle\left(27V^2 +2U^2\right)h_2, h_2\rangle+O(\epsilon^3).
\end{eqnarray*}
Similarly,
\begin{eqnarray*}
& & I(U+\epsilon h_1,V+\epsilon h_2)=  I(U,V)+2 \epsilon\left(\dpr{-\Delta U+\al U}{h_1} + 
\dpr{-\Delta V+\be V}{h_2} \right) \\
&+& \epsilon^2\left( \dpr{(-\De+\al) h_1}{h_1}+\dpr{(-\De +\be)h_2}{h_2} \right) +O(\eps^3).
\end{eqnarray*}
We now use the fact that the function 
$$
g(\eps):=\f{J(U+\epsilon h_1,V+\epsilon h_2) }{I^2(U+\epsilon h_1,V+\epsilon h_2)}
$$
has a maximum at $\eps=0$, whence $g'(\eps)=0, g''(\eps)\leq 0$. 
Hence, from the first variation, it follows that $(U, V)$ satisfies  the following system  in a weak sense  
$$\begin{cases} 
   		-\De U+\alpha U -C(\alpha,\beta)\left[(\frac{1}{9}U^2+2V^2)U+\frac{1}{3}U^2 V\right]=0 \\
   		-\De V+\beta V -C(\alpha,\beta)\left[(9V^2+2U^2)V+\frac{1}{9}U^3\right]=0.
   	\end{cases}
   $$
  where 
   $$
   C(\alpha,\beta)= \f{1}{4 J_{\max}}>0.
   $$
Focusing on the second variations we have 
\begin{eqnarray*}
& & \langle(-\Delta+\alpha-C(\alpha,\beta)(\frac{U^2}{3}+2V^2+\f{2}{3} U V))h_1,h_1\rangle+\langle(-\Delta+\beta-C(\alpha,\beta)(2 U^2+27V^2))h_2,h_2\rangle\\
& & -C(\alpha,\beta)\dpr{(8 U V+\frac{2}{3}U^2)h_1}{h_2} \geq 
 -\frac{1}{2}\left( \dpr{(-\De+\al)U}{h_1}+\dpr{-\De+\be)V}{h_2}\right)^2.
\end{eqnarray*}
 This last inequality can be rewritten as 
 $$
 \dpr{L_1 h_1}{h_1}+\dpr{L_3 h_2}{h_2}\geq 
 -\frac{1}{2}\left( \dpr{(-\De+\al)U}{h_1}+\dpr{-\De+\be)V}{h_2}\right)^2=:-\frac{1}{2} \dpr{T}{h}^2.
 $$
 where  $T=\begin{bmatrix}
 (-\De+\al)U \\
(-\De+\be) V
\end{bmatrix}.$	
With the inner product defined above we have that $\langle\cl_+h,h\rangle\geq 0.$ for $h: h\perp T.$
\end{proof}
   Based on Proposition \ref{prop:10}, we  can formulate  an  existence result  for the waves $(P,Q)$, so that they satisfy \eqref{20}. 
   \begin{proposition}
   	\label{prop:22} 
   	Assume that $\om+1>0, \mu+3\si \om>0$. Then, there exists a pair of bell-shaped functions $(P,Q)$, so that they satisfy \eqref{20}. 
   \end{proposition}
\begin{proof}
	Take $\al=\om+1>0, \be=\mu+3\si \om>0$, and apply Proposition \ref{prop:10} to this pair. This constructs $U,V$, bell-shaped and with all the extra information, mentioned in the statement and the proof of Proposition \ref{prop:10}. Finally, take 
	\begin{equation}
		\label{312} 
		P(x):=\sqrt{C(\al, \be)} U(x), Q(x):=\sqrt{C(\al, \be)} V(x).
	\end{equation}
\end{proof}
\subsection{Normalized waves}
In this section, we prove Theorem \ref{theo:40}.  To this end, for every $\la>0$, we consider the constrained variational problem 
\begin{equation}
	\label{400} 
	\begin{cases}
		E[u,v]= \frac{1}{2}\int_{\mathbb{R}} (|\nabla u|^2+
		|\nabla v|^2+|u|^2+\mu|v|^2)dx-\int_{\rone} F(u,v) dxdx\to \min \\
		M(u,v)=\int_{\mathbb{R}}  (|u|^2+3\sigma|v|^2)dx=\lambda. 
	\end{cases}
\end{equation}
where for convenience, we have introduced 
\begin{equation}
\label{001}
    F(u,v):=\frac{1}{36}|u|^4+\frac{9}{4}|v|^4+|u|^2|v|^2+\frac{1}{9}
\Re(\bar{u}^3 v).
\end{equation}

Since 
$$
|\int_{\mathbb{R}} F(u,v)dx|\leq C\|(u,v)\|_{L^4 \times L^4}^4\leq C\|(u,v)\|_{H^{\frac{1}{4}}\times H^{\frac{1}{4}}}^4\leq  C\|(u, v)\|_{H^1\times H^1} \|(u,v)\|_{L^2\times L^2}^3 \leq \epsilon \|(\nabla u,\nabla v)\|_{L^2\times L^2}^2+C_{\epsilon,\lambda},
$$
for every $\eps>0$, whence by choosing $\eps=\f{1}{2}$, we obtain 
$$
\inf_{M(u,v)=\la} E(u,v)\geq - C_{\f{1}{2},\lambda}.
$$
Thus, the constrained variational problem \eqref{400} is bounded from below and well-posed. Next, rearrangements clearly improve \eqref{400} by the Szeg\"o inequality and 
$$
\int_{\rone} |u|^2 |v|^2 \leq \int_{\rone} (u^*)^2 (v^*)^2, \Re \int_{\rone} \bar{u}^3 vdx \leq \int_{\rone} (u^*)^3 v^*dx
$$
whence $\int F(u,v) dx \leq F(u^*, v^*) dx$. Thus, we might assume, without loss of generality that the problem \eqref{400} is posed over bell-shaped functions. In particular, and due to the constraint $M(u,v)=\la$, 
$$
|u(x)|\leq C \la |x|^{-1}, |v(x)|\leq C \la |x|^{-1}.
$$
It is now easy, by Sobolev embedding and the Riesz-Kolmogrov compactness criteria to get the $L^2$ compactness of a minimizing sequence for \eqref{400}, whence the (strong $L^2$) limit of such a sequence, say $(P,Q)$ is a solution of \eqref{400}. 

Next, we derive the Euler-Lagrange equation. To this end, introduce a pair of test functions $h=(h_1, h_2)$, and consider a scalar function, 
$$
g(\eps):= E\left(\f{\sqrt{\la} (P+\eps h_1)}{\Xi(\eps)}, \f{\sqrt{\la} (Q+\eps h_2)}{\Xi(\eps)}\right)
$$
where 
$
\Xi(\eps)=\sqrt{\|P+\eps h_1\|^2+3\si \|Q+\eps h_2\|^2}.
$
Taking for simplicity $h\perp \left(\begin{matrix}
	P \\ 3 \si Q
\end{matrix}\right)$ and expanding in powers of $\eps$, we obtain 
\begin{eqnarray*}
& & 	\f{\sqrt{\la}}{\Xi(\eps)} =  1-\f{\eps^2}{2\la} \left(\int h_1^2+3\si h_2^2
		\right) +O(\eps^3);  \\
	& &  \int_{\mathbb{R}} (|\nabla (P+\eps h_1)|^2+
		|\nabla (Q+\mu \eps h_2)|^2+|(P+\eps h_1)|^2+\mu |Q+\eps h_2|^2)dx =   \\
& & 	\int_{\mathbb{R}} (|\nabla P|^2+
	|\nabla |Q|^2+P^2+\mu Q^2)dx	+ 2\eps (\dpr{-\De P +P}{h_1}+\dpr{-\De Q+\mu Q}{h_2}) +\\
	& &+ \eps^2 (\|\nabla h_1\|^2+\|h_1\|^2+\|\nabla h_2\|^2+\mu \|h_2\|^2); \\ 
	& & 	F(P+\epsilon h_1,Q+\epsilon h_2) =  F(P,Q)+
		\epsilon \left(h_1\left(\frac{1}{9}P^3+2P Q^2+\frac{1}{3}P^2Q\right)+  h_2\left(9Q^3+2P^2Q+\frac{P^3}{9}\right)\right)+\\
		& + & \epsilon^2\left(h_1\left(\frac{P^2}{6}+Q^2+\frac{1}{3}PQ\right)+\left(4PQ+\frac{1}{3}P^2\right)h_1h_2+\left(\frac{27}{2}Q^2+P^2\right)h_2^2\right)+O(\epsilon^3).
\end{eqnarray*}
Putting all this together implies the expansion 
\begin{eqnarray*}
	& & g(\eps) = g(0)+\\
	& & \eps \left(\dpr{-\De P +P-\left(\frac{1}{9}P^3+2P Q^2+\frac{1}{3}P^2Q\right)}{h_1}+\dpr{-\De Q+\mu Q-\left(9Q^3+2P^2Q+\frac{P^3}{9}\right)}{h_2}\right)+O(\eps^2).
\end{eqnarray*}
It follows that, since $(P,Q)$ is a minimizer, then $g'(0)=0$, whence 
$$
\dpr{-\De P +P-\left(\frac{1}{9}P^3+2P Q^2+\frac{1}{3}P^2Q\right)}{h_1}+\dpr{-\De Q+\mu Q-\left(9Q^3+2P^2Q+\frac{P^3}{9}\right)}{h_2} =0, 
$$
for all $h\perp \left(\begin{matrix}
	P \\ 3 \si Q
\end{matrix}\right)$. Thus, there is $\om=\om_\la$, so that 
\begin{equation}
	\label{202}   	
	\begin{cases} 
		-\De P+(\om_\lambda+1)P -\left[\frac{1}{9}P^3+2Q^2 P+\frac{1}{3}P^2Q\right]=0 \\
		-\De Q+(\mu+3\si \om_\lambda) Q -\left[9Q^3+2P^2Q+\frac{1}{9}P^3\right]=0,
	\end{cases}
\end{equation}
By Pohozaev's identities, 
$$
\omega_\lambda=\frac{-\int_{\mathbb{R}} (|\nabla P|^2+
	|\nabla Q|^2+P^2+\mu Q^2)dx+\int_{\mathbb{R}}  \frac{1}{9}P^4+9Q^4+4P^2Q^2+\frac{4}{9}
	P^3 Q)dx}{\lambda}>0.
$$
In addition, as $(P,Q)$ is a minimizer, it should be that $g''(0)\geq 0$.
\begin{eqnarray*}
& & \f{g''(0)}{2}=\langle(-\Delta+1-(\frac{P^2}{3}+2Q^2+\f{2}{3} P Q))h_1,h_1\rangle+\langle(-\Delta+\mu-(2 P^2+27Q^2))h_2,h_2\rangle \\
& &  -\dpr{(8 PQ+\frac{2}{3}P^2)h_1}{h_2}+ \omega_\lambda\epsilon^2(\int_{\mathbb{R}}h_1^2+3\sigma\int_{\mathbb{R}}h_2^2))+O(\epsilon^3)=\dpr{\cl_+h}{h}
\end{eqnarray*}
It follows that $\dpr{\cl_+h}{h}\geq 0$, whenever $h\perp \left(\begin{matrix}
	P \\ 3 \si Q
\end{matrix}\right)$. 
Finally, by the results of Proposition \ref{prop:25}, we have that for the normalized waves just constructed, the property \eqref{200} holds due to the construction. As a result, the waves $(e^{i \om t} P, e^{3 i \om t} Q)$ are spectrally stable, and the proof of Theorem \ref{theo:40} is complete. 

\section{Instability of the waves in $n=2,3$}
\label{sec:3}
In this section, we establish the instability of the waves in $n=2,3$. We need some further preparatory material, which also applies in the case $n=1$ as well,  to start our analysis with.  
   \subsection{Spectral Properties of $\cl_{\pm}$}
   
   \begin{proposition}
   	 \label{spec}
   The operators $L_{\pm}$, defined in \eqref{cpl} enjoy the following: 
   \begin{itemize}
   \item The continouous spectrum of $L_{\pm}$ is $[\min{(\alpha,\beta)},\infty).$
   \item $L_+$ has exactly one negative eigenvalue, and zero eigenvalue of multiplicity of at least $n$,  $L_+\left[\begin{matrix}
   	\p_j U \\ \p_j V
   \end{matrix}\right]=0,j=1,2,\dots,n.$
   \item $L_-\geq 0,$ with $L_-\left[\begin{matrix}
   	U \\ 3 V
   \end{matrix}\right]=0$ and moreover 
$$
\dpr{L_-  h}{h}\geq \de \|h\|^2, \forall h\perp \left[\begin{matrix}
   	U \\ 3 V
   \end{matrix}\right]
$$
   \end{itemize}
   \end{proposition}
   
   \begin{proof}
   The proof of the continuous spectrum follows directly from Weyl’s theorem. We have already establish from Proposition \ref{prop:10} that $\cl_+$ has at most one negative eigenvalue, what is left is to show the existence of one. To that end taking into account that $(U,V)$ are bell-shaped we have
   $$\langle\cl_+ {\left[\begin{matrix}
   		U \\ 3 V
   	\end{matrix}\right]},{\left[\begin{matrix}
   	U \\ 3 V
   \end{matrix}\right]}\rangle=
-C(\alpha,\beta)\int_{\rn}\left(\frac{2U^4}{9}+7U^2 V^2+18V^4+\frac{2U^3 V}{3}\right)dx<0,
$$
   thus $n(\cl_+)=1.$ Taking the spatial derivatives of \eqref{E-L}, we have $\cl_+\left[\begin{matrix}
   	\p_j U \\ \p_j V
   \end{matrix}\right]=0,j=1,2,\dots,n.$     

For the statement regarding $\cl_-,$ we start by noting that by inspection we have  $\cl_-\left[\begin{matrix}
	U \\ 3 V
\end{matrix}\right]=0$.  Now, we invoke the Perron - Frobenius property of $\cl_-$, according to which the the smallest eigenvalue  is simple and the corresponding eigenstate is a vector with non-negative functions. Thus, as $ \left[\begin{matrix}
U \\ 3 V
\end{matrix}\right]$ is such an eigenstate, $0$ is the smallest eigenvalue, which is simple, hence the claim about $\cl_-$ follows. 

Regarding the Perron-Frobenius property for $\cl_-$, this holds for any  matrix Schr\"odinger operator of the form 
$
\left(\begin{matrix}
	-\De+\al & 0 \\ 0 & -\De+\be
\end{matrix}\right)+\vec{V}
$
where $\vec{V}$ is a real-valued, bounded symmetric potential, so  in particular it holds for $\cl_-$. 
   \end{proof}
We now easily translate the results of Proposition \ref{spec} to the operators $\cl_\pm$, inroduced in \eqref{414}. 
  \begin{proposition}
	\label{spec1}
	The operators $\cl_{\pm}$, defined in \eqref{cpl} enjoy the following properties:
	\begin{itemize}
		\item 
		$$
		\si_{cont.}(\cl_{\pm})= [\min{(\alpha,\beta)},\infty).
		$$
		\item $\cl_+$ has exactly one negative eigenvalue, and zero eigenvalue of multiplicity of at least $n$,  
		$$
		\cl_+\left[\begin{matrix}
			\p_j P \\ \p_j Q
		\end{matrix}\right]=0,j=1,2,\dots,n.
	$$
		\item $\cl_-\geq 0$, 
		 with $\cl_-\left[\begin{matrix}
			P \\ 3 Q
		\end{matrix}\right]=0$ and moreover 
		$$
		\dpr{\cl_-  h}{h}\geq \de \|h\|^2, \forall h\perp \left[\begin{matrix}
			P \\ 3 Q
		\end{matrix}\right]
		$$
	\end{itemize}
\end{proposition}
   \subsection{Stability vs. instability for the eigenvalue problem \eqref{Eigen}}
   We now derive (a  necessary and sufficient) condition for stability/ instability for \eqref{Eigen}. Specifically, and based on the spectral properties of $\cl_\pm$ established in Proposition \ref{spec1}, we have the following proposition. 
   \begin{proposition}
   	\label{prop:25} 
 	If the waves $(e^{i \om t} P, e^{3 i \om t} Q)$ are spectrally stable\footnote{ i.e. the eigenvalue problem \eqref{Eigen} does not have a non-trivial solution $(\la, v): \Re\la>0, v\neq 0$}, then  $\cl_+$ is non-negative on the co-dimension one subspace $\left(\begin{matrix}
   		P \\ 3\si Q
   	\end{matrix}\right)^\perp$, i.e. 
   	\begin{equation}
   		\label{200} 
   	\dpr{\cl_+ h}{h}\geq 0, \forall h\perp \left(\begin{matrix}
   		P \\ 3\si Q
   	\end{matrix}\right). 
   	\end{equation}
   Equivalently, spectral instability occurs if and only if there exists $h\perp \left(\begin{matrix}
   	P \\ 3\si Q
   \end{matrix}\right)$, so that $\dpr{\cl_+ h}{h}<0$. 
   \end{proposition}
\begin{proof}
	We first need a technical step to eliminate the problem posed by $\f{1}{\si}$ in the definition of $\cj$. To this end, denote for simplicity 
	$
	I_\si=\left(\begin{matrix}
		1 & 0 \\ 0 & \f{1}{\si}
	\end{matrix}\right)
	$, so that 
	$$
	\cj=\left(\begin{matrix}
		{\mathbf 0}_2 & I_\si \\ - I_\si & 	{\mathbf 0}_2 
	\end{matrix}\right).
	$$
	We can rewrite the eigenvalue problem \eqref{Eigen} as follows 
	\begin{equation}
		\label{205} 
	\begin{cases}
		I_\si \cl_- v_2 = \la v_1 \\
		 I_\si \cl_+ v_1 = -\la v_2.
	\end{cases}
	\end{equation}
Introducing $v_j\to \sqrt{I_\si} v_j, j=1,2$, we can further equivalently rewrite  
\eqref{205} as 
\begin{equation}
	\label{206} 
	\begin{cases}
		\sqrt{I_\si}  \cl_- \sqrt{I_\si} v_2 = \la v_1 \\
		\sqrt{I_\si} \cl_+ \sqrt{I_\si} v_1 = -\la v_2.
	\end{cases}
\end{equation}
Introducing the new operators $\tilde{\cl}_\pm:=\sqrt{I_\si}  \cl_\pm \sqrt{I_\si} $, we arrive at the more standard NLS form 
\begin{equation}
	\label{207} 
	\begin{cases}
	\tilde{\cl}_- v_2 = \la v_1 \\
	\tilde{\cl}_+ v_1 = -\la v_2
	\end{cases}
\end{equation}
Note that since $\cl_-$ is positive on the subspace $\left(\begin{matrix}
	P \\ 3 Q
\end{matrix}\right)^\perp$, it is easy to see that $\tilde{\cl}_-\left(\begin{matrix}
P \\ 3\sqrt{\si}  Q
\end{matrix}\right)=0$ and $\tilde{\cl}_-$  is positive on the subspace $\left(\begin{matrix}
P \\ 3\sqrt{\si}  Q
\end{matrix}\right)^\perp$. 

Applying now $\tilde{\cl}_-$ on the second equation of \eqref{207}, we arrive at the eigenvalue problem 
\begin{equation}
	\label{209} 
	\tilde{\cl}_- \tilde{\cl}_+ v_1=-\la^2 v_1,
\end{equation}
which is clearly equivalent to the original eigenvalue problem \eqref{Eigen} in the sense that instability occurs exactly when \eqref{209} has solution $(\la, v_1): \Re \la>0, v_1\in (H^4(\rn))^2\neq 0$. Assume that $\la\neq 0, v_1\neq 0$. 

Taking dot product of \eqref{209} with $Z_\si:=\left(\begin{matrix}
	P \\ 3\sqrt{\si}  Q
\end{matrix}\right)$ implies that $v_1\perp Z_\si$. Since $\tilde{\cl}_-$ is positive on $Z_\si^\perp$, we have that $\tilde{\cl}_-: Z_\si^\perp\to Z_\si^\perp$ is well-defined and strictly positive, whence we can introduce 
$$
h: \sqrt{\tilde{\cl}_-}  h=v_1.
$$
As we can view the eigenvalue problem \eqref{209} on the subspace $Z_\si^\perp$, the latest allows us to rewrite it as 
\begin{equation}
	\label{210} 
	\sqrt{\tilde{\cl}_-} \tilde{\cl}_+ \sqrt{\tilde{\cl}_-} h = -\la^2 h.
\end{equation}
From \eqref{210}, it is immediately clear that $-\la^2$ is real, as eigenvalue of the  self-adjoint operator $\sqrt{\tilde{\cl}_-} \tilde{\cl}_+ \sqrt{\tilde{\cl}_-}$, so $\la\in {\mathbf R}\cup i {\mathbf R}$. In particular, instabilities, if any, present only as positive eigenvalues. 

Assume spectral instability, that is \eqref{210} has a solution, which must be $\la: \la>0$. Taking dot product with $h$ implies 
$$
0>-\la^2 \dpr{h}{h}=\dpr{	\sqrt{\tilde{\cl}_-} \tilde{\cl}_+ \sqrt{\tilde{\cl}_-} h}{h}=\dpr{\tilde{\cl}_+ v_1}{v_1}.
$$
Thus, we have produced a vector (namely $v_1=\sqrt{\tilde{\cl}_-} h \perp Z_\si$), so that $\dpr{\tilde{\cl}_+ v_1}{v_1}<0$. But this means 
$$
\dpr{\cl_+ \sqrt{I}_\si v_1}{\sqrt{I}_\si  v_1}<0,
$$
where $\sqrt{I}_\si v_1 \perp \left(\begin{matrix}
	P \\ 3 \si Q
\end{matrix}\right)$. One direction is established. 

For the opposite direction, assume that there is $h: h\perp \left(\begin{matrix}
	P \\ 3 \si Q
\end{matrix}\right), \dpr{\cl_+ h}{h}<0$. This implies that for 
$h=:\sqrt{I_\si} h_\si$, we have $\perp Z_\si$ and 
$$
\dpr{\tilde{\cl}_+ h_\si}{h_\si}<0.
$$
Then, introduce $v: \sqrt{\tilde{\cl}_-} v=h_\si$. It follows that $\dpr{\sqrt{\tilde{\cl}_-} \tilde{\cl}_+ \sqrt{\tilde{\cl}_-} v}{v}<0$, whence the self-adjoint operator  $\sqrt{\tilde{\cl}_-} \tilde{\cl}_+ \sqrt{\tilde{\cl}_-}$ has a negative eigenvalue. As a consequence,   \eqref{210} has solutions and spectral instability follows. 
\end{proof}
\subsection{ Instability of the waves in $n>1$: Proof of Theorem \ref{theo:20}} 
We have the following proposition. 
\begin{proposition}
	\label{prop:45} 
	Let  $P,Q$ are classical solutions of \eqref{20}. Introduce 
	$$
	H(x)=\left(\begin{matrix}
		H_1 \\
		H_2
	\end{matrix}\right):=\left(\begin{matrix}
		\f{n}{2} P+x\cdot \nabla P \\
		\f{n}{2}  Q+ x\cdot \nabla Q
	\end{matrix}\right).
$$
 Then, $H\perp  \left(\begin{matrix}
  P \\
 3 Q
 \end{matrix}\right)$ and 
$$
\begin{cases}
\dpr{\cl_+ H}{H}<0, n=3, \\  \dpr{\cl_+ H}{H}=0, n=2.
\end{cases}
$$
\end{proposition}
As a consequence of Proposition \ref{prop:25}, we immediately conclude spectral instability in the case $n=3$, while an additional argument is needed to establish the instability in the case $n=2$. 
\begin{proof}
	Apply $x\cdot \nabla $ to the elliptic equation \eqref{20}. Using the commutation identity 
	\begin{equation}
		\label{235} 
			-x\cdot(\De)= -\De(x\cdot) +2 \De,
	\end{equation}
	 we obtain 
	 \begin{equation}
	 	\label{215} 
	 	\cl_+ \left(\begin{matrix}
	 	x\cdot \nabla P \\
	 	 x\cdot \nabla Q
	 	\end{matrix}\right)=-2 \left(\begin{matrix}
	 	\De  P \\
	 	\De Q
 	\end{matrix}\right).
	 \end{equation}
	A direct calculation yields 
	\begin{equation}
		\label{220} 
		\cl_+ \left(\begin{matrix}
			 P \\
			 Q
		\end{matrix}\right)=-\left(\begin{matrix}
		2(-\De P+(\om+1) P)\\
		2(-\De Q+(\mu+3\si \om) Q)
		\end{matrix}\right).
	\end{equation}
	Alltogether, 
	\begin{equation}
		\label{225} 
		\cl_+ H= \left(\begin{matrix}
			(n-2)\De P-n(\om+1) P)\\
			(n-2) \De Q - n(\mu+3\si \om) Q
		\end{matrix}\right).
	\end{equation}
	Taking dot product with $H$ (noting  that by construction $H\perp \left(\begin{matrix}
		 P\\
		0
	\end{matrix}\right), \left(\begin{matrix}
	0\\
	Q
\end{matrix}\right)$), we have 
$$
\dpr{\cl_+ H}{H}=(n-2) \left(\dpr{\De P}{x\cdot P}+ \dpr{\De Q}{x\cdot Q}\right).
$$
Clearly, for the case $n=2$, we have $\dpr{\cl_+ H}{H}=0$, while for the case $n=3$, we use the commutation identity \eqref{235} to integrate by parts. We obtain 
$$
\dpr{\De P}{x\cdot P} = -\int_{\rthree} |\nabla P|^2, \dpr{\De Q}{x\cdot Q} = -\int_{\rthree} |\nabla Q|^2.
$$
	It follows that for $n=3$, $\dpr{\cl_+ H}{H}=-\left(\int_{\rthree} |\nabla P|^2+|\nabla Q|^2\right)<0$. 
\end{proof}
We finally are ready to state the spectral  instability result. 
\begin{proposition}
	\label{prop:59} 
	Let $n=2,3$ and  $P,Q$ are classical solutions of \eqref{20}, which enjoy the spectral properties listed in Proposition \ref{spec1}. Then, 
	\begin{itemize}
		\item For $n=3$,  the pair $(e^{i \om t} P, e^{3 i \om t} Q)$  is spectrally unstable. 
		\item  For $n=2$, and under the extra assumption $\mu\neq 3\si$,  $(e^{i \om t} P, e^{3 i \om t} Q)$ is unstable as well. 
		\item 
		Finally, in the special case $n=2, \mu=3\si$,   the waves  $(e^{i \om t} P, e^{3 i \om t} Q)$ are  spectrally stable. 
	\end{itemize}
\end{proposition}
\begin{proof}
	As we have alluded above, by combining  the results of Proposition \ref{prop:25} and Proposition \ref{prop:45}, we have spectral instability in the case $n=3$. 
	
	For the case $n=2$ and $\mu\neq 3\si$, assume spectral stability, for a contradiction. By Proposition \ref{prop:25} (see \eqref{200}), it must be that $\cl_+$ is non-negative on the subspace $\left(\begin{matrix}
		P \\ 3Q
	\end{matrix}\right)^\perp$. By the proof of Proposition \ref{prop:25}, this is in turn equivalent to the non-negativity of the operator\footnote{Note that this operator acts invariantly on $\left(\begin{matrix}
		P \\ 3\sqrt{\si} Q
	\end{matrix}\right)^\perp$, due to the presence of  $\sqrt{\tilde{\cl}_-}$} $\sqrt{\tilde{\cl}_-} \tilde{\cl}_+ \sqrt{\tilde{\cl}_-}$ on the subspace $\left(\begin{matrix}
	P \\ 3\sqrt{\si}  Q
\end{matrix}\right)^\perp$. 

That is, with the projection $\Pi_{\left(\begin{matrix}
	P \\ 3\sqrt{\si}  Q
\end{matrix}\right)^\perp}: L^2 \to \left(\begin{matrix}
P \\ 3\sqrt{\si}  Q
\end{matrix}\right)^\perp$, we have 
$$
\sqrt{\tilde{\cl}_-} \Pi_{\left(\begin{matrix}
	P \\ 3\sqrt{\si}  Q
\end{matrix}\right)^\perp}\tilde{\cl}_+ \Pi_{\left(\begin{matrix}
P \\ 3\sqrt{\si}  Q
\end{matrix}\right)^\perp}\ \sqrt{\tilde{\cl}_-} :  \left(\begin{matrix}
P \\ 3\sqrt{\si}  Q
\end{matrix}\right)^\perp\to \left(\begin{matrix}
P \\ 3\sqrt{\si}  Q
\end{matrix}\right)^\perp
$$
 is non-negative. On the other hand, 
$$
\inf_{Z\perp \left(\begin{matrix}
		P \\ 3\sqrt{\si}  Q
	\end{matrix}\right): \|Z\|=1} \dpr{ \sqrt{\tilde{\cl}_-} \tilde{\cl}_+ \sqrt{\tilde{\cl}_-} Z}{Z}=0.
$$
In fact, one can take $Z_0=\|\tilde{\cl}_-^{-\f{1}{2}} H\|^{-1} \tilde{\cl}_-^{-\f{1}{2}} H$, with $H$ as in Proposition \ref{prop:45}. It must be then, that $0$ is an eigenvalue for the self-adjoint operator $\sqrt{\tilde{\cl}_-} \tilde{\cl}_+ \sqrt{\tilde{\cl}_-} $, with eigenvector $Z_0$. But this means, as $ \sqrt{\tilde{\cl}_-}$ is invertible on $\left(\begin{matrix}
	P \\ 3\sqrt{\si}  Q
\end{matrix}\right)^\perp$, that 
$$
\Pi_{\left(\begin{matrix}
		P \\ 3\sqrt{\si}  Q
	\end{matrix}\right)^\perp} \tilde{\cl}_+ \sqrt{\tilde{\cl}_-} Z_0=0,
$$
which is the same as 
$$
\Pi_{\left(\begin{matrix}
		P \\ 3 \si Q
	\end{matrix}\right)^\perp} \cl_+ H=0.
$$
In view of \eqref{225} however, 
$$
0=\Pi_{\left(\begin{matrix}
		P \\ 3 \si Q
	\end{matrix}\right)^\perp}   \cl_+ H =-3\Pi_{\left(\begin{matrix}
	P \\ 3 \si Q
\end{matrix}\right)^\perp}  \left(\begin{matrix}
(\om+1) P)\\
(\mu+3\si \om) Q
\end{matrix}\right)=-3(\mu-3\si) \Pi_{\left(\begin{matrix}
	P \\ 3 \si Q
\end{matrix}\right)^\perp} \left(\begin{matrix}
0 \\ Q
\end{matrix}\right)\neq 0,
$$
as long as $P\neq 0$ and $\mu\neq 3\si$. 
The contradiction implies that the wave $(e^{i \om t} P, e^{3 i \om t} Q)$ is spectrally unstable for $n=2$ as well. 

For the case $n=2$, $\mu=3\si$,  since $n(\cl_+)=1$ and $\cl_-\geq 0$, the stability is decided by the Vakhitov-Kolokolov  quantity $\dpr{\cl_+^{-1} \left(\begin{matrix}
		P \\ 3 \si Q
	\end{matrix}\right)}{\left(\begin{matrix}
	P \\ 3 \si Q
\end{matrix}\right)}$, see the eigenvalue problem \eqref{207}. 

Our calculations (see \eqref{225} with $n=2, \mu=3\si$), showed that we have 
$$
\cl_+ H=-2(\om+1) \left(\begin{matrix}
	P \\ 3\si Q
\end{matrix}\right). 
$$
It follows that 
$$
\dpr{\cl_+^{-1} \left(\begin{matrix}
		P \\ 3 \si Q
	\end{matrix}\right)}{\left(\begin{matrix}
		P \\ 3 \si Q
	\end{matrix}\right)}=-\f{1}{2(\om+1)} \dpr{H}{\left(\begin{matrix}
P \\ 3 \si Q
\end{matrix}\right)}=0.
$$
This shows (marginal) spectral stability as in the case of 2D cubic NLS, where there is an extra pair of (generalized) eigenvectors at zero. Indeed, the  algebraic multiplicity of the zero eigenvalue for \eqref{Eigen} in this case is eight, instead of the usual six. 
\end{proof}
\section{Strong instability}
\label{sec:4} 
In this section, for a fixed $\om>\max(-1, -\f{\mu}{3\si})$, we consider the waves $P,Q$ constructed in Proposition \ref{prop:22}, specifically see \eqref{312}. In particular, they are solutions of the Weinstein minimization problem \eqref{42}. 

We shall also need to restrict our parameters as follows $\si=3, \mu=3\si=9$. With these specific values, we record the energy and mass as follows 
\beqn
E(u,v) &=& \frac{1}{2}\int_{\rn} (|\nabla u|^2+
|\nabla v|^2+|u|^2+9 |v|^2)dx-\int_{\rn}  \frac{1}{36}|u|^4+\frac{9}{4}|v|^4+|u|^2|v|^2+\frac{1}{9}
\Re(\bar{u}^3 v)dx\\
M(u,v) &=& \int_{\rn}  (|u|^2+9 |v|^2)dx.
\eeqn

We start with the following proposition, which provides the necessary virial identities, which are ultimately the basis for our instability by blow up results. 
\begin{proposition}
	\label{viral}
	let $n=2,3 ,$ $(u,v)\in H^1(\mathbb{R}^n)\cap L^2(\mathbb{R}^n,|x|^2dx) $ and $\sigma=3,\mu=3\si=9$ and let  
	$$
	V(t)=\int |x|^2\left(|u(t)|^2+9|v(t)|^2\right) dx 
	$$
	 Then 
	\begin{equation}
		\label{400} 
			V''(t)=\int _{\rone^n}\left(8|\nabla u|^2+8|\nabla v|^2-\frac{2n}{9}|u|^4-18n|v|^4-8n|v|^2|u|^2-\frac{8n}{9}
		\Re(\bar{u}^3 v) \right) dx
		=:16R(u,v). 
	\end{equation}
	Where $$R(u,v):=\frac{1}{2}\left(\int _{\rone^n}|\nabla u|^2+|\nabla v|^2-\frac{n}{36}|u|^4-\frac{9n}{4}|v|^4-n|v|^2|u|^2-\frac{n}{9}
	\Re(\bar{u}^3 v) dx\right).
	$$
	
\end{proposition}
\begin{proof}
	This follows from Proposition 4.4 in \cite{OP}.
\end{proof}
Next, introduce the $L^2-$scale invariant transformation,   
$u^\la(x):=(\la^{\frac{n}{2}}u(\la x)$.  We compute, 
\begin{eqnarray*}
	\partial_\la E(u^\la,v^\la) &=& \la\int_{\rone^n}|\nabla u|^2+|\nabla v|^2-n\la^{n-1}\int_{\rone^n}F(u,v)dx, \\ 
	\partial_\la ^2E(u^\la,v^\la) &=& \int_{\rone^n}|\nabla u|^2+|\nabla v|^2-n(n-1)\la^{n-2}\int_{\rone^n}F(u,v)dx.
\end{eqnarray*}
Clearly, 
\begin{eqnarray}
	\label{300} 
& & 	\frac{\la}{2}\partial_\la E(u^\la,v^\la)=R(u^\la,v^\la), \\
	\label{002}
& & 	\|\nabla P\|_{L^2(\rone^n)}^2+ \|\nabla Q\|_{L^2(\rone^n)}^2=n\int_{\rone^n}F(P,Q)dx.
\end{eqnarray}
where \eqref{002}  is a consequence of \eqref{Poh1}, \eqref{Poh2}, \eqref{Poh3}. 

Following the approach of Ohta, \cite{O}, we introduce the sets 
\begin{eqnarray*}
	\mathcal{A} &=& \left\{(u,v) : E(u,v)<E(P,Q),\int P ^2+9Q^2 =\int u^2+9u^2,\int F(P,Q) <\int F(u,v) \right\},\\
	\mathcal{B} &=& \{(u,v)\in \mathcal{A}:  R(u,v)<0\}.
\end{eqnarray*}
We now split our considerations in the cases $n=2$ and $n=3$. We start with the harder case $n=3$. 
\subsection{The case $n=3$}
\begin{lemma}
	Let $n=3$. Then $(P^\la,Q^\la)\in \mathcal{B},$ for all $\la>1.$
\end{lemma}
\begin{proof}
	By the invariance of the scaling, we have $\int _{\rone^3}P ^2+9Q^2 dx=\int _{\rone^3}(P^\la)^2+9(Q^\la)^2$. 
	Next,  observe that as $F$ is a quartic term, 
	$$
	 \int_{\rone^3}F(P^\la,Q^\la)dx=\la^3\int_{\rone^3}F(P,Q)dx
	 >\int_{\rone^3}F(P,Q)dx
	 $$
	  for $\la>1$.  
We have that 
\begin{eqnarray*}
	2 R(P^\la, Q^\la) &=& \la^2 \int _{\rthree}|\nabla P|^2+|\nabla Q|^2-
\la^3 	\int_{\rthree} \frac{1}{12}|P|^4+\frac{27}{4}|Q|^4+3 |P|^2|Q|^2+\frac{1}{3}
	P^3 Q =\\
	&=& \la^2 \int _{\rthree}|\nabla P|^2+|\nabla Q|^2 - 3\la^3  \int F(P, Q) \\
	&=& 3\la^2(1-\la) \int _{\rthree}|\nabla P|^2+|\nabla Q|^2<0, 
\end{eqnarray*}
where we have used \eqref{002} and  $\la>1$.  By the formula \eqref{300}, we conclude that 
$$
\frac{\la}{2}\partial_\la E(P^\la,Q^\la)=R(P^\la,Q^\la)<0,
$$
 whence the scalar function $\la\to E(P^\la,Q^\la)$ is decreasing for $\la>1$. In particular, 
 $$
 E(P^\la, Q^\la)<E(P,Q),
 $$
whenever $\la>1$. We have thus shown that $(P^\la, Q^\la)\in \ca$ for $\la>1$. Recall that we have also shown that $R(P^\la, Q^\la)<0$, whence $(P^\la, Q^\la)\in \cb$ as well. 
\end{proof}
Next, we shall need the following technical lemma. 
\begin{lemma}
	\label{003}
	If $(u,v)\in (H^1(\rone^3)\times H^1(\rone^3))$ such that 
	\begin{equation}
		\label{340} 
			\int _{\rone^3}P ^2+9Q^2 dx=\int _{\rone^3}u^2+9u^2dx;  \ \ \int_{\rone^3}F(u,v)dx=\int_{\rone^3}F(P,Q)dx.
	\end{equation}
	 Then $E(P,Q)\leq E(u,v)$. 
\end{lemma}
\begin{proof}
	Recall that $(P,Q)$ solve the minimzation problem \eqref{42}, for\footnote{Note the formula \eqref{312} which relates $(P,Q)$ to the minimizers $(U,V)$. As the problem \eqref{42} is homogeneously invariant, so $(P,Q)$ are also solutions}  under  $\al=\om+1$ and $\be=9(1+\om)$. 
	Moreover, by the homogeneity of \eqref{42}, $(P,Q)$ solves the following constrained minimization problem 
	\begin{equation}
		\label{320} 
		\begin{cases}
			K[u,v]:=\f{1}{2} \int_{\rthree} |\nabla u|^2 + |\nabla v|^2 + (1+\om)|u|^2+9(1+\om)|v|^2 \to \min \\
			 \int_{\rthree} F(u,v)=\int_{\rthree} F(P,Q)
		\end{cases}
	\end{equation}
	Let now $(u,v)$ satisfy the constraints \eqref{340}. As $M(u,v)=M(P,Q)$ and 
	$ \int_{\rthree} F(u,v)=\int_{\rthree} F(P,Q)$, it follows from the fact that $(P,Q)$ solves \eqref{320} that 
	$$
E(P,Q)+\int_{\rthree} F(P,Q)+\f{\om}{2} M(P,Q)	=K(P,Q)\leq K(u,v)= E(u,v)+\int_{\rthree} F(u,v)+\f{\om}{2} M(u,v).
	$$
	Cancelling out the equal terms leads to the desired inequality  $E(P,Q)\leq E(u,v)$.

\end{proof}
Next, we show that the set $\ca$ is invariant under the flow. 
\begin{lemma}
	\label{03}
	The set $\mathcal{A} $ is invariant under the flow of \eqref{10}, at least until a potential blow-up. 
\end{lemma}
\begin{proof}
	Take initial data $(u_0, v_0)\in \ca$. That is $E(u_0, v_0)<E(P,Q), M(u_0, v_0)=M(P,Q), \int_{\rthree} F(P, Q) <\int F(u_0, v_0)$. Evolving this initial data in time, by using the conservation of mass and energy, we see that 
	\begin{eqnarray}
		\nonumber
	& & 	M(u(t), v(t))=M(u_0, v_0)=M(P,Q); \\ 
\label{380}	
	& & 	E(u(t), v(t))=E(u_0, v_0)<E(P,Q); 
	\end{eqnarray}
In addition, the continuity of the solution map guarantees that for some small time interval, $\int_{\rthree} F(P, Q) <\int_{\rthree} F(u(t), v(t))$, say in $t\in (0, \de)$. 

	Hence, the only possibility for $(u(t), v(t))\notin \ca$ for some $t_0>0$ is that eventually, 
	$$
	\int _{\rthree}F(u(t_0), v(t_0))=\int_{\rthree} F(P, Q).
	$$
	 But then, we are in a position to apply Lemma \ref{003}. This implies $E(u(t_0), v(t_0))\geq E(P,Q)$, in contradiction with \eqref{380}. Thus, $(u(t), v(t))\in \ca$ for as long as the potential blow up occurs. 
\end{proof}
Note that for  $n=3$,
\begin{equation}
	\label{004}
	E(u,v)-R(u,v)=\frac{1}{2}M(u,v)+\frac{1}{2}\int_{\rone^3}F(u,v)dx.
\end{equation}
We now need another technical lemma. 
\begin{lemma}
	\label{00003}
	Let $n=3$ and if $(u,v)\in (H^1(\rone^3)\times H^1(\rone^3))$ such that $R(u,v)\leq 0, M(u,v)=M(P,Q),$ $ \int_{\rone^3}F(P,Q)dx<\int_{\rone^3}F(u,v)dx.$ Then 
	$$
	E(P,Q)\leq E(u,v)-R(u,v).
	$$
\end{lemma}

\begin{proof}
	Using the identity \eqref{004}, the conservation of mass and the assumptions  on $(u,v)$, we have
	\begin{eqnarray*}
		E(P,Q) &=& \frac{1}{2}M(P,Q)+\frac{1}{2}\int_{\rone^3}F(P,Q)dx=\frac{1}{2}M(u,v)+\frac{1}{2}\int_{\rone^3}F(P,Q)dx\leq \\
		&\leq &  \frac{1}{2}M(u,v)+\frac{1}{2}\int_{\rone^3}F(u,v)=E(u,v)-R(u,v)
	\end{eqnarray*}
\end{proof}
We are then ready to show that the set $\cb$ is also  invariant under the flow. 
\begin{lemma}
	\label{005}
	Assume $n=3,$ then $\mathcal{B}_\om$ is invariant under the flow of \eqref{10}.
\end{lemma}

\begin{proof}
	Let $(u_0, v_0)\in \cb$, so in particular $(u_0, v_0)\in \ca$. 
	
	 Since $\mathcal{A}$ is invariant we need to show $R(u,v)<0$. 
	Suppose there exist $t_1\in(0,T_{max})$ such that $R(u(t_1),v(t_1))\geq0.$ Then by continuity of $t\to R(u(t),v(t))$ there is $t_0\in(0,t_1]$ such that $R(u(t_0),v(t_0))=0.$ Then by Lemma \ref{03}, 
	$$
	M(u(t_0),v(t_0))=M(P,Q); \ \ \int_{\rone^3}F(P,Q)dx<\int_{\rone^3}F(u(t_0),v(t_0))dx.
	$$
	 Thus by Lemma \ref{00003} we have
	$$
	E(P,Q)\leq E(u(t_0),v(t_0))-R(u(t_0),v(t_0))=E(u(t_0),v(t_0).
	$$
	
	This is a contradiction to the statement in Lemma \ref{003} that $E(u(t_0),v(t_0))<E(P,Q)$, hence $\cb$ is invariant as stated. 
\end{proof}
We are now ready for the proof of the strong instability of the waves $(P,Q)$, that is Theorem \ref{21} in the case $n=3$. 

To this end, let $(u_0,v_0)\in \mathcal{B}_\om$ and $(u(t),v(t))$ be the solution \eqref{10} with initial data $(u_0,v_0).$ Then by Lemma \ref{005} $(u(t),v(t))\in \mathcal{B}_\om$ for all $t\in[0,T_{max}).$ Also the virial identity Proposition \ref{viral} and Lemma \ref{003} together with conservation of energy implies
$$\frac{1}{16}V''(t)=R(u(t),v(t))\leq E(u(t),v(t))-E(P,Q)=E(u_0,v_0)-E(P,Q)<0$$
for every $t\in [0,T_{max}).$ This implies that $T_{max}<\infty$. 

\subsection{The case $n=2$}
We have the following simplified version of the virial identity. 
\begin{lemma}
	For $n=2,\si=3, \mu=9$, take  initial data $(u_0,v_0)\in (H^1(\rone^2)\times H^1(\rone^2))$. Then, for all times, up to a potential blow-up time, 
	 we have 
	\begin{equation}
		\label{410} 
		V''(t)=16E(u(t),v(t))-8M(u(t),v(t))=16E(u_0,v_0)-8M(u_0,v_0).
	\end{equation}
\end{lemma}
\begin{proof}
	The formula  \eqref{410} is just a particular instance of \eqref{400} in the case $n=2$, where we have used the Pohozhaev's identities. Finally, we used the conservation of mass and energy to conclude that this is actually a constant in time. 
\end{proof}
Note that in this case, it would suffice to find initial data $u_0, v_0$ so that 
\begin{equation}
	\label{430} 
	16E(u_0,v_0)-8M(u_0,v_0)<0.
\end{equation}
Another observation is that for $u_0=P, v_0=Q$, we have that $V(t)=const$, whence $V''(t)=0$, so we get\footnote{this can of course also be obtained from the Pohozhaev's identities}  $16 E(P,Q)-8 M(P, Q)=0$. Cleary, in order to produce data with the property \eqref{430} is to perturb the wave $(P,Q)$. To this end, take a pair of real test functions $h_1, h_2$ and consider the perturbation $(P+\epsilon h_1,Q+\epsilon h_2)$. Expanding up to first  order in\footnote{Observe that the zero order in $\eps$ vanishes due to the relation $16 E(P,Q)=8 M(P,Q)$}  $\eps$ we obtain 
\begin{eqnarray*}
& & 	E(P+\epsilon h_1,Q+\epsilon h_2)-\f{1}{2} M(P+\epsilon h_1,Q+\epsilon h_2)=\\
	&=& 
	\eps \left(\dpr{-\De P - (\f{1}{9}P^3+2 P Q^2+\f{1}{3} P^2 Q)}{h_1}\right)+ 
	\eps\left(\dpr{-\De Q - (9 Q^3+2 P^2 Q+\f{1}{9} P^3)}{h_2}\right)+O(\eps^2)\\
	&=& -(\om+1) \eps \dpr{\left(\begin{array}{c} P \\ 9 Q \end{array}\right)}{\left(\begin{array}{c} h_1 \\ h_2 \end{array}\right)} +O(\eps^2).
\end{eqnarray*}
 Selecting,  say $h_1=P, h_2=Q$,  and taking into account that $\om>-1$, implies that for all small $0<\eps<<1$, we have 
 $E(P + \epsilon P,Q  + \epsilon Q )<\f{1}{2} M(P+\epsilon P,Q+\epsilon Q)$. Thus, solutions with this initial data blow up in finite time.

 Theorem \ref{21} is thus proved in full.


\end{document}